\documentclass[reqno0,12pt]{amsart}
\pdfoutput=1

\usepackage{dsfont}
\usepackage{bbm}
\usepackage[nodate]{datetime}
\usepackage[hmargin=24mm,top=28mm,bottom=28mm,a4paper]{geometry}
\usepackage{color}
\usepackage{subfig}
\usepackage{fancyhdr}
\usepackage{amsmath,amssymb,amsfonts,amsthm}
\usepackage{amstext}
\usepackage{amsmath}
\usepackage{amssymb}
\usepackage{enumerate}
\usepackage{amsbsy}
\usepackage{amsopn}
\usepackage{bbm,amsthm}
\usepackage{amscd}
\usepackage{amsxtra}
\usepackage{todonotes}
\usepackage{soul}

\newtheorem{theorem}{Theorem}[section]
\newtheorem*{theorem*}{Theorem}
\newtheorem{lemma}{Lemma}[section]

\theoremstyle{remark}

\newcommand{\ind}{\mathds{1}}

\begin{document}
\title{Sign changes of the partial sums of a random multiplicative function II}
\author{Marco Aymone}
\begin{abstract}
We study two models of random multiplicative functions: Rademacher random multiplicative functions supported on the squarefree integers $f$, and Rademacher random completely multiplicative functions $f^*$. We prove that the partial sums $\sum_{n\leq x}f^*(n)$ and $\sum_{n\leq x}\frac{f(n)}{\sqrt{n}}$ change sign infinitely often as $x\to\infty$, almost surely. The case $\sum_{n\leq x}\frac{f^*(n)}{\sqrt{n}}$ is left as an open question and we stress the possibility of only a finite number of sign changes, with positive probability. 
\end{abstract}
\maketitle
\section{Introduction}
\subsection{Main results and motivation} A Rademacher random multiplicative function $f$ is defined as follows: Over the primes, the values $f(p)$ are \textit{i.i.d.} random variables taking $\pm 1$ values with probability $1/2$ in each instance (Rademacher random variables), and at the other values of $n$,
$$f(n):=\mu^2(n)\prod_{p|n}f(p),$$
where $\mu$ is the M\"obius function.

By a Rademacher random completely multiplicative function $f^*$, we mean that
$$f^*(n)=\prod_{p^a\| n}f(p)^a.$$

In a previous paper, Aymone, Heap and Zhao \cite{aymone_sign_changes} proved that the partial sums $\sum_{n\leq x}f(n)$ of a Rademacher random multiplicative function change sign infinitely often\footnote{Here we say that a function $M(x)$ changes sign infinitely often as $x\to\infty$ if at a certain point $x_0$, $M(x)$ is not always non-negative or non-positive for all $x\geq x_0$.} as $x\to\infty$, almost surely. The key idea was to show that, almost surely
\begin{align*}
\lim_{\sigma\to1/2^+}\int_1^\infty\left(\sum_{n\leq x}f(n)\right)\frac{dx}{x^{1+\sigma}}=0,\\
\lim_{\sigma\to1/2^+}\int_1^\infty\left|\sum_{n\leq x}f(n)\right|\frac{dx}{x^{1+\sigma}}=\infty,
\end{align*}  
and these visibly capture an infinite number of sign changes.

Our task here is to study sign changes for weighted sums of $f$ and $f^*$:
$$\sum_{n\leq x}\frac{f(n)}{n^\alpha},\,\sum_{n\leq x}\frac{f^{*}(n)}{n^\alpha}, \;\alpha\geq 0.$$
In both cases, we have that for $\alpha>1/2$, the partial sums converge to a non-vanishing Euler product, and hence, it can be shown that these partial sums become positive for all $x$ sufficiently large, almost surely. Therefore, we restrict ourselves to the cases $0\leq\alpha\leq 1/2$.

Our motivation comes from the study of the sign changes of the partial sums of the deterministic counterparts to $f$ and $f^*$: the  M\"obius $\mu$ and the Liouville $\lambda$. For $0\leq \alpha <1/2$, the partial sums $\sum_{n\leq x}\mu(n)n^{-\alpha}$ change sign infinitely often as $x\to\infty$, due to the analytic properties of $1/\zeta(s)$ and to the Landau's oscillation theorem. Answering a question of P\'olya, Haselgrove \cite{haselgrove} proved that $\sum_{n\leq x}\lambda(n)$ changes sign infinitely often as $x\to\infty$. In the range $0<\alpha<1/2$, the case of an infinite number of sign changes of $\sum_{n\leq x}\lambda(n)n^{-\alpha}$ recasts questions about the non-trivial Riemann zeros, and at the present moment, unconditional results are unknown, see \cite{trudgian_sign_changes} by Mossinghoff and Trudgian. In a paper by Humphries \cite{humphries_sign_changes}, it can be deduced from his results that under RH and other hypothesis on Riemann zeros, the partial sums $\sum_{n\leq x}\lambda(n)n^{-\alpha}$ change sign infinitely often as $x\to\infty$. 

In both cases $\sum_{n\leq x}\mu(n)n^{-1/2}$ and $\sum_{n\leq x}\lambda(n)n^{-1/2}$ the problem of an infinitude of sign changes seems to be intricate and even assuming standard hypothesis, conditional results are unknown. In \cite{trudgian_sign_changes} is left as an open problem to prove or disprove that $\sum_{n\leq x}\lambda(n)n^{-1/2}$ is negative for all $x\geq 17$. A remarkable consequence of this fact is that this implies that all non-trivial zeros of the Riemann zeta function have multiplicity at most $2$, see Theorem 2.4 of \cite{trudgian_sign_changes}. In \cite{humphries_sign_changes} it is showed that under standard assumptions on the non-trivial zeros of the Riemann zeta function, the set of $x$ for which $\sum_{n\leq x}\frac{\lambda(n)}{\sqrt{n}}$ is negative has total logarithmic asymptotic density.

Taking all of these facts into account, for the random counterparts $f$ and $f^*$ we were able to proof the following results.

\begin{theorem}\label{teorema principal 1} Let $f$ be a Rademacher random multiplicative function. Then for each $0\leq\alpha\leq 1/2$, $\sum_{n\leq x}\frac{f(n)}{n^\alpha}$ changes sign infinitely often as $x\to\infty$, almost surely.
\end{theorem} 
\begin{theorem}\label{teorema principal 2}
Let $f^*$ be a Rademacher random completely multiplicative function. Then for each $0\leq\alpha< 1/2$ we have that $\sum_{n\leq x}\frac{f^*(n)}{n^\alpha}$ changes sign infinitely often as $x\to\infty$, almost surely.
\end{theorem}

The case $\sum_{n\leq x}\frac{f^*(n)}{\sqrt{n}}$ is left as an open question, and it is possible that these sums change sign only a finite number of times with positive probability, here we give an heuristic reason. In a recent paper, Angelo and Xu \cite{xu_positivity} (and later improvements by Kerr and Klurman \cite{klurman_positivity}) studied the weighted sums $\sum_{n\leq x}\frac{f^*(n)}{n}$ in the context of a ``\textit{random Turán conjecture}''. They proved that the probability of $\sum_{n\leq x}\frac{f^*(n)}{n}$ being positive for all $x\geq1$ is at least $1-10^{-45}$. This implies a certain positivity for harmonic sums of $f^*$ and it may be possible that such result may be extended for $\sum_{n\leq x}\frac{f^*(n)}{n^\sigma}$, where $\sigma>1/2$. Further, after a convolution argument, $f^*=f\ast \ind_{PS}$, where $\ind_{PS}(n)$ is the indicator that $n$ is a perfect square. Therefore, 
$$\sum_{n\leq x}\frac{f^*(n)}{\sqrt{n}}\approx \sum_{n\leq x}\frac{f(n)}{\sqrt{n}}\log(x/n).$$
The later is known as a Riesz mean, and at some instances these averages stabilizes oscillatory properties of an oscillating arithmetic function providing a finite number of sign changes for these averages. 

 \subsection{Proof method} Our proof follows closely the lines of \cite{aymone_sign_changes}. The difference here is that, in contrast with unweighted sums of a Rademacher random multiplicative function, we have the almost sure limits
\begin{align*}
\lim_{\sigma\to1/2^+}\int_1^\infty\left(\sum_{n\leq x}\frac{f(n)}{\sqrt{n}}\right)\frac{dx}{x^{1/2+\sigma}}=\infty,\\
\lim_{\sigma\to1/2^+}\int_1^\infty\left|\sum_{n\leq x}\frac{f(n)}{\sqrt{n}}\right|\frac{dx}{x^{1/2+\sigma}}=\infty,
\end{align*} 
and this alone does not capture an infinite number of sign changes, unless we make the two divergences above quantitative, and that these divergences differ in size. This can, indeed, be achieved by using the following quantitative statement by Harper:
\begin{theorem}[Harper \cite{harpergaussian}, page 25]\label{teorema Harper} Let $p$ run over the primes and $(f(p))_p$ be a sequence of independent Rademacher random variables. For fixed $A>3$, almost surely, there exists a sequence $\sigma_k\to1/2^+$ such that
\begin{align*}\label{equacao harper}
\sup_{1\leq t\leq 2\left(\log\left(\frac{1}{\sigma_k-1/2}\right)\right)^2}\left(\sum_{p}\frac{f(p)\cos(t\log p)}{p^{\sigma_k}}-2\log\log\left(\frac{1}{\sigma_k-1/2}\right)\right) \\ 
\geq\log\left(\frac{1}{\sigma_k-1/2}\right)-A\log\log\left(\frac{1}{\sigma_k-1/2}\right).
\end{align*} 
\end{theorem}

\subsection{Background} Here we are intended to give a list of results in the literature (perhaps non-complete) on random multiplicative functions. The first result on partial sums of random multiplicative functions is due to Wintner \cite{wintner}, where he considered a question by Levy \cite{levy} and showed that $\sum_{n\leq x}f(n)$ is almost surely $\Omega(x^{1/2-\epsilon})$ and $\ll x^{1/2+\epsilon}$, for any $\epsilon>0$.

Later, these results have been improved by Erd\H{o}s \cite{erdosuns} and subsequently by Hal\'asz \cite{halasz}. Later Basquin \cite{basquinn} and Lau, Tenenbaum and Wu \cite{tenenbaum2013} made a significant improvement by proving that, almost surely, 
$$\sum_{n\leq x}f(n)\ll \sqrt{x}(\log\log x)^{2+\epsilon},$$
for any $\epsilon>0$.

In this year of 2023, Caich \cite{rachidupperbound} proved that we can replace $(\log \log x)^{2+\epsilon}$ above by $(\log \log x)^{1/4+\epsilon}$. This combined with Harper's Omega bound \cite{harper_largefluctations} gives a sharp description of the size of the fluctuations of the partial sums of $f$: 
$$\sum_{n\leq x}f(n)=\Omega(\sqrt{x}(\log\log x)^{1/4-\epsilon}),$$ 
for any $\epsilon>0$, almost surely. 

Central limit Theorems also have been studied in various settings: \cite{chatterjeemult}, \cite{harpermult}, \cite{houghmult}, \cite{klurman_rand_chowla}, \cite{xu_clt}.

The problem of moments also have been extensively studied, with the remarkable result of Harper that shows that the first moment exhibits better than squareroot cancellation \cite{harperfirstmoment}. Other related works (including related models):
 \cite{aymone_JLMS}, \cite{benatar_rodgers}, \cite{seip_helsons}, \cite{hardy_weighted} \cite{harper_high_moments}, \cite{harper_ashkhan}, \cite{heap_lindqvist}, \cite{xu_firstmoment}.
 
And another viewpoints of investigation \cite{xu_positivity}, \cite{aymone_sign_changes}, \cite{aymonebiased}, \cite{kalmynin_quadratic}, \cite{klurman_positivity}, \cite{najnudel}.

\section{Preliminaries}
\subsection{Notation} Here $p$ stands for a generic prime number. We use the standard Vinogradov notation $f(x)\ll g(x)$ or Landau's $f(x)=O(g(x))$ whenever there exists a constant $c>0$ such that $|f(x)|\leq c|g(x)|$, for all $x$ in a set of parameters. When not specified, this set of parameters is an infinite interval $(a,\infty)$ for sufficiently large $a>0$. The standard $f(x)=o(g(x))$ means that $f(x)/g(x)\to0$ when $x\to a$, where $a$ could be a complex number or $\pm \infty$. For $g(x)>0$ for all $x$, we say that $f(x)=\Omega(g(x))$ if $\limsup_{x\to\infty}\frac{|f(x)|}{g(x)}>0$. 
\subsection{Some Lemmas}
By the pioneering work of Wintner \cite{wintner}, we have that $F(s):=\sum_{n=1}^\infty\frac{f(n)}{n^s}$ converges for all $Re(s)>1/2$, almost surely. His proof can be modified to show that the same is true for $F^*(s)=\sum_{n=1}^\infty\frac{f^*(n)}{n^s}$. Moreover we have:
\begin{lemma}\label{lemma Euler product} For all $Re(s)>1/2$, almost surely we have
\begin{align*}
F(s)=&\prod_{p}\left(1+\frac{f(p)}{p^s}\right),\\
F^*(s)=&\prod_{p}\left(1-\frac{f(p)}{p^s}\right)^{-1}.
\end{align*}
\end{lemma}

This leads to the almost sure identities
\begin{equation}\label{equacao exponencial F}
F(s)=\exp\left(\sum_{p}\frac{f(p)}{p^s}-\frac{1}{2}\log\zeta(2s)+O(1)\right),
\end{equation}
\begin{equation}\label{equacao exponencial F*}
F^*(s)=\exp\left(\sum_{p}\frac{f(p)}{p^s}+\frac{1}{2}\log\zeta(2s)+O(1)\right),
\end{equation}
where the $O(1)$ term above is actually an analytic function that is uniformly bounded in $Re(s)\geq 1/3+\epsilon$, for small fixed $\epsilon>0$. For a proof of these results we refer reader to Lemma 2.4 of \cite{aymone_sign_changes} and the references therein.
To conclude this section, we have
\begin{lemma}\label{lemma upper bounding RDS}
Let $(f(p))_p$ be independent Rademacher random variables. Then, for all $\epsilon>0$, we have almost surely as $\sigma\to1/2^+$
$$\sum_p\frac{f(p)}{p^\sigma}\ll \left(\log\left(\frac{1}{\sigma-1/2}\right)\right)^{1/2+\epsilon}.$$
\end{lemma}
For a proof we refer reader to Lemma 2.1 and Remark 2.1 of \cite{aymone_sign_changes}. 
\section{Proof of the main results}
\subsection{The case $\sum_{n\leq x}\frac{f(n)}{\sqrt{n}}$} By the upper bound for the partial sums of a Rademacher random multiplicative function \cite{basquinn} and \cite{tenenbaum2013},
$$\sum_{n\leq x}f(n)\ll \sqrt{x}(\log\log x)^{2+\epsilon},$$
we have that by partial summation
$$\sum_{n\leq x}\frac{f(n)}{\sqrt{n}}\ll (\log x)^2, \;a.s.$$
Further refinements of this last upper bound could, perhaps, be achieved, by following the lines of \cite{aymone_JLMS} and more recently of \cite{hardy_weighted}. 

Now, for $Re(s)=\sigma>1/2$,
$$F(s):=\sum_{n=1}^\infty \frac{f(n)}{n^s}=\sum_{n=1}^\infty \frac{f(n)}{\sqrt{n}}\frac{1}{n^{s-1/2}}=(s-1/2)\int_{1}^\infty \left(\sum_{n\leq x}\frac{f(n)}{\sqrt{n}} \right)\frac{dx}{x^{s+1/2}}.$$

By equation \eqref{equacao exponencial F} and Lemma \ref{lemma upper bounding RDS} combined with the estimate $\log\zeta(2\sigma)=\log(1/(2\sigma-1))+O(1)$ valid for $\sigma\to 1/2^+$, we obtain that as $\sigma\to 1/2^+$, almost surely we have $F(\sigma)= (2\sigma-1)^{1/2+o(1)}$.
Hence, we almost surely have as $\sigma\to1/2^+$
\begin{equation}\label{equacao integral f rad weighted}
\int_{1}^\infty \left(\sum_{n\leq x}\frac{f(n)}{\sqrt{n}} \right)\frac{dx}{x^{\sigma+1/2}}=\frac{1}{(2\sigma-1)^{1/2+o(1)}}.
\end{equation} 

On the other hand, for $t\geq 1$
$$\left|\frac{F(\sigma+it)}{t}\right|\ll \int_{1}^\infty \left|\sum_{n\leq x}\frac{f(n)}{\sqrt{n}} \right|\frac{dx}{x^{\sigma+1/2}}.$$

Now, by equation \eqref{equacao exponencial F},
$$|F(s)|=\exp\left(\sum_{p}\frac{f(p)\cos(t\log p)}{p^\sigma}-\frac{1}{2}Re(\log\zeta(2s))+O(1)\right).$$

We recall that $|\log \zeta(\sigma+it)|\leq \log\log|t|+O(1)$, for $\sigma\geq 1$ and $t\geq 2$. By Theorem \ref{teorema Harper}, almost surely we have an infinite sequence of points $\sigma_k\to1/2^+$ and $t_k\geq 1$ such that
$$\left|\frac{F(\sigma_k+it_k)}{t_k}\right|\geq\exp\left((1+o(1))\log\left(\frac{1}{\sigma_k-1/2}\right)\right)\gg \frac{1}{(\sigma_k-1/2)^{0.9}}.$$
Therefore, almost surely, there exists a sequence $\sigma_k\to1/2^+$ such that
\begin{equation}\label{equacao divergencia modulo rad weighted}
\int_{1}^\infty \left|\sum_{n\leq x}\frac{f(n)}{\sqrt{n}} \right|\frac{dx}{x^{\sigma_k+1/2}}\gg \frac{1}{(\sigma_k-1/2)^{0.9}}.
\end{equation}

To complete the argument, we see that the quantities at the left of \eqref{equacao integral f rad weighted} and \eqref{equacao divergencia modulo rad weighted} have different size, and so the partial sums $\sum_{n\leq x}\frac{f(n)}{\sqrt{n}}$ cannot be always non-negative or non-positive, there always must be infinite sign changes, almost surely.

\subsection{The case $\sum_{n\leq x}\frac{f^*(n)}{n^\alpha}$} Let then $0\leq\alpha<1/2$. Let $F^*(s)$ be the Dirichlet series of $f^*$. Similarly as in the previous case, by equation \eqref{equacao exponencial F*} and Lemma \ref{lemma upper bounding RDS}, we have that
$$F^*(\sigma)=\frac{1}{(2\sigma-1)^{1/2+o(1)}},$$
almost surely as $\sigma\to1/2^+$.

By partial summation and integration,
$$F^*(s)=(s-\alpha)\int_1^\infty \left(\sum_{n\leq x}\frac{f^*(n)}{n^\alpha}\right)\frac{dx}{x^{s+1-\alpha}}.$$
Therefore, almost surely as $\sigma\to 1/2^+$
\begin{equation}\label{equacao rad compl diverg sem mod}
\int_1^\infty \left(\sum_{n\leq x}\frac{f^*(n)}{n^\alpha}\right)\frac{dx}{x^{\sigma+1-\alpha}}=\frac{1+o(1)}{(1/2-\alpha)}\frac{1}{(2\sigma-1)^{1/2+o(1)}}.
\end{equation}
And similarly to \eqref{equacao divergencia modulo rad weighted}, almost surely there exists a sequence $\sigma_k\to1/2^+$ such that
\begin{equation}\label{equacao rad compl diverg mod}
\int_1^\infty \left|\sum_{n\leq x}\frac{f^*(n)}{n^\alpha}\right|\frac{dx}{x^{\sigma+1-\alpha}}\gg\frac{1}{(2\sigma_k-1)^{0.99}}.
\end{equation}
Just as before, \eqref{equacao rad compl diverg sem mod} and \eqref{equacao rad compl diverg mod} captures infinite sign changes, almost surely.

\subsection{The case $\sum_{n\leq x}\frac{f(n)}{n^\alpha}$} This case can be treated similarly as the previous one, although we can proceed in a more elementary way as in \cite{aymone_sign_changes}, that is, Theorem \ref{teorema Harper} is not needed here.

The proof of Theorems \ref{teorema principal 1} and \ref{teorema principal 2} are now complete.

\noindent\textbf{Acknowledgements} This project started while I was a visiting Professor at Aix-Marseille Université, and I am warmly thankful for their kind hospitality, and to CNPq for supporting this visit -- grant PDE no. 400010/2022-4 (200121/2022-7). This project also was funded by CNPq grant Universal no. 403037/2021-2 and by FAPEMIG grant Universal no. APQ-00256-23.

{\small{\sc \noindent
Departamento de Matem\'atica, Universidade Federal de Minas Gerais (UFMG), Av. Ant\^onio Carlos, 6627, CEP 31270-901, Belo Horizonte, MG, Brazil.} \\
\textit{Email address:} \verb|aymone.marco@gmail.com| }

\end{document}